\pdfoutput=1
\RequirePackage{ifpdf}
\ifpdf 
\documentclass[pdftex]{sigma}
\else
\documentclass{sigma}
\fi

\usepackage{mathtools}
\usepackage[shortlabels]{enumitem}

\usepackage{tikz}
\usetikzlibrary{arrows,calc,decorations.pathreplacing,decorations.markings,intersections,shapes.geometric,through,fit,shapes.symbols,positioning,decorations.pathmorphing}

\numberwithin{equation}{section}

\newtheorem{Theorem}{Theorem}[section]
\newtheorem*{Theorem*}{Theorem}
\newtheorem{Corollary}[Theorem]{Corollary}
\newtheorem{Lemma}[Theorem]{Lemma}
\newtheorem{Proposition}[Theorem]{Proposition}

\theoremstyle{definition}
\newtheorem{Definition}[Theorem]{Definition}

\newtheorem{Example}[Theorem]{Example}
\newtheorem{Remark}[Theorem]{Remark}
\newtheorem{remarks}[Theorem]{Remarks}

\newcommand\bA{{\mathbb A}}

\newcommand\bD{{\mathbb D}}

\newcommand\bF{{\mathbb F}}
\newcommand\bG{{\mathbb G}}
\newcommand\bH{{\mathbb H}}

\newcommand\bK{{\mathbb K}}

\newcommand\bN{{\mathbb N}}
\newcommand\bO{{\mathbb O}}
\newcommand\bP{{\mathbb P}}
\newcommand\bQ{{\mathbb Q}}
\newcommand\bR{{\mathbb R}}
\newcommand\bS{{\mathbb S}}
\newcommand\bT{{\mathbb T}}
\newcommand\bU{{\mathbb U}}

\newcommand\bZ{{\mathbb Z}}

\newcommand\cA{{\mathcal A}}
\newcommand\cB{{\mathcal B}}

\newcommand\cF{{\mathcal F}}

\newcommand{\cat}[1]{\textsc{#1}}

\DeclareMathOperator{\id}{id}

\DeclareMathOperator{\Aut}{\mathrm{Aut}}

\newcommand{\xrightarrowdbl}[2][]{%
 \xrightarrow[#1]{#2}\mathrel{\mkern-14mu}\rightarrow
}

\begin{document}

\allowdisplaybreaks

\newcommand{\arXivNumber}{2206.11997}

\renewcommand{\PaperNumber}{073}

\FirstPageHeading

\ShortArticleName{Prescribed Graphon Symmetries and Flavors of Rigidity}

\ArticleName{Prescribed Graphon Symmetries\\ and Flavors of Rigidity}

\Author{Alexandru CHIRVASITU}

\AuthorNameForHeading{A.~Chirvasitu}

\Address{Department of Mathematics, University at Buffalo, Buffalo, NY 14260-2900, USA}
\Email{\mail{achirvas@buffalo.edu}}

\ArticleDates{Received November 06, 2025, in final form August 03, 2026; Published online August 10, 2026}

\Abstract{We prove that an arbitrary compact metrizable group can be realized as the automorphism group of a graphon; this is a continuous analogue to Frucht's theorem recovering arbitrary finite groups as automorphism groups of finite graphs. The paper also contains a~number of results on the persistence of transitivity of a compact-group action upon passing to a limit of graphons. Call a compact group $\mathbb{G}$ graphon-rigid if, whenever it acts transitively on each member $\Gamma_n$ of a convergent sequence of graphons, it also acts transitively on the limit $\lim_n \Gamma_n$. We show that for a compact Lie group~$\mathbb{G}$ graphon rigidity is equivalent to the identity component $\mathbb{G}_0$ being semisimple; as a partial converse to a~result of Lov\'{a}sz and Szegedy, this is also equivalent to weak randomness: the property that the group have only finitely many irreducible representations in each dimension. Similarly, call a compact group~$\mathbb{G}$ image-rigid if for every compact Lie group $\mathbb{H}$ the images of morphisms $\mathbb{G}\to \mathbb{H}$ form a~closed set (of closed subgroups, in the natural topology). We prove that graphon rigidity implies image rigidity for compact groups that are either connected or profinite, and the two conditions are equivalent (and also equivalent to being torsion) for profinite abelian groups.}


\Classification{22C05; 28A33; 60B05; 22F50; 20B27; 54E35; 54E45; 54E50; 54E70}

\section{Introduction}
Graphons are continuous and/or probabilistic analogues of graphs, intended to capture the latter's long-term, limiting behavior. One of several slight variations (e.g., \cite[Section~2.1]{ls}; we review standard spaces briefly in Section~\ref{se:prel}).

\begin{Definition}
 A {\it graphon} $\Gamma=(X,W,\mu)$ consists of a standard probability space $(X,\mu)$ and a~measurable symmetric function $W\colon X\times X\to [0,1]$.
\end{Definition}

Aspects of the rich, multi-faceted ensuing theory are discussed in \cite{ls-limdense,ls-reg,ls,sz-limker} (and their references) as well as the excellent~\cite{lvs-bk}. The topic of interest here is that of symmetries (or automorphisms) of a graphon. The automorphism group $\Aut(\Gamma)$ is introduced in passing in \cite[Section~13.5]{lvs-bk}, and discussed more extensively in~\cite{ls}. Most importantly for the discussion below, it is always compact \cite[Theorem 10]{ls}.

One direction the discussion below takes in the general spirit of an ``inverse'' symmetry problem: given a group (enjoying various properties, as appropriate), realize it as the automorphism group of a structure (here, graphons). The pattern recurs frequently:
\begin{itemize}\itemsep=0pt
\item Every finite group is the automorphism group of some finite graph by a classical result of Frucht's \cite{frucht}.
\item That result can be improved by imposing various constraints on the graphs realizing the given group (connectivity, chromatic number, etc.) \cite[Theorem 1.2]{sabid-props}.
\item And in fact {\it arbitrary} (possibly infinite) groups can be realized as graph automorphism groups \cite[Introduction, Theorem]{sab-inf}.
\item Adjacently to graphs, (convex) {\it polytopes} also realize arbitrary finite groups, see \cite[Theorem 1.1]{doign}, \cite[Theorems 1 and 2]{sw}, \cite[Theorem 3.1]{ssw}.
\item There are variations on this theme: the group might have a distinguished involution which is required to act as the central symmetry in a {\it symmetric} convex polytope; again, arbitrary data of this nature can be realized \cite[Theorem 2.1]{cls}.
\item And then there is of course the celebrated {\it inverse Galois problem}, of realizing finite groups as Galois groups of finite field extensions of $\bQ$ (still unsolved in full generality) \cite{ser-gal}.
\end{itemize}

The graphon analogue of all of this is Theorem \ref{th:prscr} below.

\begin{Theorem*}
 Every compact metrizable group is the automorphism group of some graphon.
\end{Theorem*}

Or in other words, since \cite[Theorem 10]{ls} says that the automorphism group of a graphon is compact (and it is obviously metrizable), that result has as strong a converse as one could possibly hope for.

Another point of interest is the limiting behavior of graphon automorphism groups (limits being understood throughout as those defined via \emph{morphism densities}, see \cite[Section~2.1]{ls} and the brief recollection opening Section~\ref{se:grphn}). The issue comes up in \cite{ls}, where \cite[Examples~34 and~35]{ls} show that if
\begin{equation*}
 \Gamma_n\xrightarrow[\quad n\quad]{}\Gamma
\end{equation*}
is a convergent sequence of graphons with $\Gamma_n$ all admitting transitive actions by a group $\bG$, the limit $\Gamma$, though again transitive under $\Aut(\Gamma)$, need not admit a transitive action by the {\it same} group $\bG$.

By contrast, \cite[Theorem~37]{ls} identifies (following \cite[Section~1.8]{sz-limker}) a class of groups for which transitivity {\it does} transport over from the individual $\Gamma_n$ to the limit $\Gamma=\lim_n \Gamma_n$: the {\it weakly random} compact groups, i.e., those having only finitely many irreducible representations in every dimension.

To simplify the awkward language, we refer to compact groups $\bG$ for which transitive actions transfer over from the members $\Gamma_n$ of a convergent sequence to the limit $\Gamma=\lim_n \Gamma_n$ as {\it graphon-rigid} (see Definition \ref{def:rig}). Weak randomness and graphon rigidity are both (as the latter's name suggests) rigidity properties of sorts. In light of the above the question naturally arises of whether (and to what extent) weak randomness is {\it necessary} for graphon rigidity. The two turn out to coincide for Lie groups, and to in fact be equivalent to yet another (much more?) familiar property in that case. To paraphrase the more verbose Theorem \ref{th:lierig} as follows.

\begin{Theorem*}
 A compact Lie group is graphon-rigid precisely when it is weakly random or, equivalently, when its identity connected component is semisimple.
\end{Theorem*}

Being Lie is essential: Example~\ref{ex:z2s} shows that it is perfectly possible for infinite {\it abelian} compact groups to be graphon-rigid. Being abelian, these are maximally far from weak randomness: {\it all} of their irreducible representations are 1-dimensional.

The various (counter)examples suggest a related notion of rigidity (Definition~\ref{def:downrig}): call a~compact group $\bG$ {\it image-rigid} if for every compact Lie group $\bH$, the images of the morphisms $\bG\to \bH$ form a closed set in the natural compact Hausdorff topology (Remark~\ref{re:haustop}) on the space $\mathcal{CLG}(\bH)$ of closed subgroups of $\bH$.

This property, in some cases, appears to be more pertinent to graphon rigidity than weak randomness. Aggregating Theorem~\ref{th:profab} and Proposition~\ref{pr:prof2rig}, we have the following.

\begin{Theorem*}
 Let $\bG$ be a compact group.
 \begin{itemize}\itemsep=0pt
 \item If $\bG$ is graphon-rigid, then it is image-rigid provided it is either connected or profinite.
 \item If furthermore $\bG$ is abelian and profinite, then graphon and image rigidity are equivalent. \qedhere
 \end{itemize}
\end{Theorem*}

In the short Section \ref{se:graphing}, we give a possible definition of the automorphism group of a {\it graphing}, in partial answer to a question asked in \cite[Section~5]{lvs-grphng}. That group will not, in general, be compact, but the definition does make it invariant under taking {\it full subgraphings} in the sense of \cite[Section~2]{lvs-grphng}.

\section{Preliminaries}\label{se:prel}

We assume some basic background on measure theory, as covered in, say, \cite{pet,bog}; the references will be more precise when needed. As is customary in the literature (e.g., \cite[Section~2.2]{bog}), for a measure $\mu$ we abbreviate {\it almost everywhere} or {$\mu$-almost everywhere} to {\it a.e.} or {\it $\mu$-a.e.}. Recall also \cite[Definition~2.1.3]{bog} that a function
$\varphi\colon   (X,\cB_X)\to (Y,\cB_Y)$
between spaces equipped with {\it $\sigma$-algebras} \cite[Definition~1.2.2]{bog} is {\it measurable} when
$\varphi^{-1}(S)\in \cB_X$, $\forall S\in \cB_Y$.

We work virtually exclusively with {\it probability} measures, and all probability spaces $(X,\cB,\mu)$ are {\it complete} (see \cite[Section~1.1]{pet} or \cite[Definition 1.5.10]{bog}): subsets of measurable sets $N\subset X$ with $\mu(N)=0$ are again measurable. The following objects feature prominently in the sequel (see, e.g., \cite[Definitions~9.2.1 and~9.4.6]{bog}).

\begin{Definition}\label{def:std}\quad
 \begin{enumerate}[(a)]\itemsep=0pt
 \item Two probability spaces $(X,\cB_X,\mu_X)$ and $(Y,\cB_Y,\mu_Y)$ are {\it isomorphic mod 0} if
 \begin{itemize}\itemsep=0pt
 \item there are sets of $M\subset X$ and $N\subset Y$ of measure zero,
 \item and a measurable, measure-preserving bijection
$
 X\setminus M\to Y\setminus N$,
 \item whose inverse is also measurable.
 \end{itemize}
 \item A {\it standard} probability space is a (complete) probability space that is isomorphic mod~0 to the unit interval with its usual Lebesgue measure appropriately scaled, together, perhaps, with a sequence of point masses.
 \end{enumerate}
\end{Definition}

\begin{Remark}\label{re:std}
 Standard probability spaces are also (essentially)
 \begin{itemize}\itemsep=0pt
 \item the {\it Lebesgue spaces} of \cite[Section~1.4, Definition 4.5]{pet};
 \item the {\it Lebesgue--Rohlin spaces} of \cite[Definition~9.4.6 and Theorem~9.4.7]{bog};
 \item the {\it Lebesgue--Rohlin spaces} of \cite[Definition~7]{haez} (see \cite[Proposition~6 and Remarks following it]{haez}).
 \end{itemize}
 The terminology also matches \cite[Appendix~A.3.1]{lvs-bk}.
\end{Remark}

\begin{Remark}\label{re:mod0}
 The concept of mod-0 morphism
$(X,\cB_X,\mu_X)\to (Y,\cB_Y,\mu_Y)$
 is potentially ambiguous, with at least three possible interpretations. It might mean
 \begin{enumerate}[(a)]\itemsep=0pt
 \item\label{item:1} a measurable, measure-preserving map defined on a full-measure set $X\setminus M$ as in Definition~\ref{def:std}, with two such maps identified if they agree a.e.;
 \item\label{item:2} a measurable, measure-preserving map defined {\it everywhere} on $X$, again identifying two such maps if they agree almost everywhere;
 \item\label{item:3} or, finally, a measure-preserving morphism
$\cB_Y\to \cB_X$
 of $\sigma$-algebras \cite[Section~1.4.C]{pet} (note the direction reversal: such a map is to be thought of as the preimage function $\varphi^{-1}$ induced by $\varphi$ as in~\ref{item:1} or~\ref{item:2}). Or, in the language of \cite[Section~1.4.C]{pet}, a morphism of {\it measure algebras}.
 \end{enumerate}
 For standard spaces the ambiguity is only apparent:
 \begin{itemize}\itemsep=0pt
 \item On the one hand the notions listed above are clearly progressively ``weaker'' as listed, in the sense that an object as in \ref{item:1} produces one as in \ref{item:2}, etc. Being standard plays no role here.
 \item On the other hand, for standard spaces a $\sigma$-algebra morphism as in~\ref{item:3} arises from~\ref{item:2} by \cite[Section~1.4, Theorem 4.7]{pet}.
 \item And the passage from~\ref{item:2} back to~\ref{item:1} is illustrated in the proof of \cite[Proposition~6]{haez}.
 \end{itemize}

 As all of our spaces are standard unless noted otherwise, we can make free use of these mod-0 concepts in whatever form is most convenient in any given situation. The discussion also applies to {\it iso}morphisms, {\it auto}morphisms, etc. For instance, for a standard probability space $(X,\cB,\mu)$ we can speak of its {\it automorphism group} $\Aut(X,\cB,\mu)$: it consists of mod-0 classes of automorphisms in the sense of the preceding discussion.
\end{Remark}

\section{Prescribed graphon symmetries}\label{se:grphn}

The many references throughout the text to graphon-sequence convergence \smash{$W_n\!\xrightarrow[\raisebox{2pt}{\scriptsize $n$}]{}\! W$} are always to the same notion, based on \emph{morphism densities} (see \cite[Section~11.1]{lvs-bk} or \cite[Section~2.1]{ls}):
\begin{equation*}
 \begin{aligned}
& (X,W_n,\mu_n)\xrightarrow[\quad n\quad]{\quad}(X,W,\mu)
  \iff
 \forall (\text{finite simple graph $F$} )
 \bigl(t(F,W_n)\xrightarrow[\quad n\quad]{\quad}t(F,W)\bigr),\\
& t(F,W)
:=
 \int_{X^V}\prod_{ij\in E}W(x_i,x_j)\prod_{i\in V}\mathrm{d}\mu(x_i),
 \end{aligned}
\end{equation*}
with $V$ and $E$ denoting the vertex and edge sets of $F$ respectively. The same convergence can be obtained via the \emph{cut distance} \cite[Section~8.2.2]{lvs-bk} on weak isomorphism classes of graphons by \cite[Lemmas~10.23 and~10.32]{lvs-bk}, but we will not employ that metric explicitly in the sequel. The same quantities $t(F,W)$ also determine the \emph{weak isomorphism} class of a graphon (by definition \cite[Section~2.1, p.~139]{ls}), so matters of convergence depend only on such isomorphism classes.

Recall from \cite[Section~2.3]{ls} (also \cite[Section~13.3]{lvs-bk}):

\begin{Definition}\label{def:nbhd}
 The {\it neighborhood distance} $r_W$ attached to a graphon $(X,W,\mu)$ is
 \begin{equation*}
 r_W(x,y):=\|W(x,-) - W(y,-)\|_1 = \int_{X}|W(x,z)-W(y,z)|\, \mathrm{d}\mu(z).
 \end{equation*}
 A graphon $(J,W,\mu)$ is {\it pure} if it is complete under the neighborhood distance and its probability measure $\mu$ has full support.
\end{Definition}

We work with metric probability spaces $(X,d,\mu)$ as graphons, as explained in \cite[Example~13.17]{lvs-bk}. These turn out to automatically be pure, as follows from the following simple remark.

\begin{Proposition}\label{pr:cpct}
 A graphon $\Gamma=(X,W,\mu)$ with
 \begin{itemize}\itemsep=0pt
 \item $X$ compact Hausdorff second-countable,
 \item $W$ continuous and
 \begin{equation}\label{eq:xtowx}
 X \ni x\mapsto W(x,-)\in C(X) :=  \text{continuous functions on $X$}
 \end{equation}
 one-to-one,
 \item and $\mu$ fully supported
 \end{itemize}
 is pure.
\end{Proposition}
\begin{proof}
 For compact Hausdorff spaces being second-countable is equivalent to metrizability \cite[Section~34, Exercise~3]{mnk}, and the continuity of $W$ then entails its {\it uniform} continuity \cite[Theorem~27.6]{mnk}. It follows that~\eqref{eq:xtowx} is continuous when equipping the right-hand side with the uniform topology, so the identity map is continuous from the original (compact Hausdorff) topology on~$X$ and the $r_W$-metric topology.

 Note also that $r_W$ is indeed a metric (rather than a pseudometric) because we are assuming that~\eqref{eq:xtowx} is injective. It follows that
$\id\colon  X\to (X,r_W)$
 is a bijective and continuous map from a~compact space to a Hausdorff space, and must thus be a homeomorphism. This then implies that the metric space $(X,r_W)$ is complete (being compact \cite[Theorem 45.1]{mnk}), and purity follows from the additional assumption that $\mu$ is fully supported.
\end{proof}

In particular, taking for the kernel $W$ the distance function $d$ of a compact metric space $(X,d)$ as in \cite[Example 13.17]{lvs-bk}, we obtain

\begin{Corollary}\label{cor:cpct}
For a compact metric space $(X,d)$ with fully-supported probability measure~$\mu$, the resulting graphon $(X,d,\mu)$ is pure.
\end{Corollary}

The following observation will be implicit in some of the discussion below.

\begin{Lemma}\label{le:actcont}
 Under the hypotheses of Proposition~{\rm \ref{pr:cpct}}, the automorphism group of $\Gamma$ acts continuously on~$X$ in the latter's original compact Hausdorff topology.
\end{Lemma}
\begin{proof}
 Indeed, it acts continuously on the metric space $(X,r_W)$, and we saw in the course of the proof that metric space is homeomorphic to the original~$X$.
\end{proof}

The following result is a graphon analogue of arbitrary groups being realizable as graph automorphism groups \cite[p.~64, Theorem]{sab-inf}. There are approaches to graphon automorphism groups ostensibly different from \cite[Definition~7]{ls}: compare Definitions~\ref{def:origgrphn} and~\ref{def:autobis} below. The ambiguities will be elided in Corollary~\ref{cor:sameauto}, which confirms that the resulting objects coincide.

\begin{Theorem}\label{th:prscr}
 A compact group is the automorphism group of a graphon if and only if it is metrizable.
\end{Theorem}
\begin{proof}
 On the one hand, \cite[Theorem 10]{ls} shows that the automorphism group of a graphon is compact, while \cite[discussion following Definition~7]{ls} explains how to metrize it.

 Conversely, \cite[Theorem~1.2]{mell} shows that an arbitrary compact metrizable group $\bG$ is of the form $\bG\cong \Aut(X,d)$ for a compact metric space $(X,d)$. Now recast $(X,d)$ as a pure graphon, as follows:
 \begin{itemize}\itemsep=0pt
 \item rescale $d$ so that it takes values in $[0,1]$;
 \item and consider a fully-supported probability measure $\mu$ on $X$, invariant under $\bG$; one can always be produced by averaging {\it any} fully-supported measure $\nu$ (extant, as $X$ is second-countable compact Hausdorff) against the Haar measure $\mu_{\bG}$ of~$\bG$:
 \begin{equation*}
 \mu(Y) :=
 \int_{\bG}\nu(gY)\, \mathrm{d}\mu_{\bG}(g).
 \end{equation*}
 \end{itemize}
 By the very choice of $\mu$ as invariant under the \emph{entirety} of $\bG\cong \Aut(X,d)$, the canonical embedding $\Aut(X,\mu,d)\lhook\joinrel\to \Aut(X,d)$ is onto and hence a compact-group isomorphism. Our graphon is thus now $(X,\mu,d)$, pure by Corollary \ref{cor:cpct}, with automorphism group precisely the original~$\bG$.
\end{proof}

\subsection{Revisiting a definition}\label{subse:bckgrphn}

The approach taken in \cite[Definition 7]{ls} to defining the automorphism group of a graphon is different from that suggested by Remark \ref{re:mod0}: rather than define the automorphisms as classes modulo zero (and perhaps modulo automorphisms that are in other ways ``trivial''), loc.\ cit.\ first normalizes the graphon to a pure one and then isolates concrete automorphisms thereof.

\begin{Definition}\label{def:origgrphn}
 Let $(X,W,\mu)$ be a pure graphon. Its {\it automorphism group} $\Aut_{\rm LS}(X,W,\mu)$ consists of those measure-preserving bijections $\varphi\colon X\to X$ such that
 \begin{equation*}
 \forall x\in X,\qquad W(\varphi x,  \varphi \cdot) = W(x,\cdot)\text{ almost everywhere}.
 \end{equation*}
\end{Definition}

The `LS' subscript (for Lov\'{a}sz--Szegedy) is temporary; we will see that the definition specializes to something more akin to \cite[Section~1]{vh}.

\begin{Definition}\label{def:autobis}
 Let $\Gamma := (X,W,\mu)$ be a graphon (pure or not). Its {\it automorphism group} $\Aut(\Gamma)$ is defined as follows:
 \begin{itemize}\itemsep=0pt
 \item Consider the group $\Aut(X,\mu)$ of mod-0 equivalence classes of measure-preserving automorphisms of the standard probability space $(X,\mu)$.
 \item Set
 \begin{equation*}
 \widetilde{\Aut}(\Gamma):=\{g\in \Aut(X,\mu)\mid  W(g\, \cdot,  g\, \cdot) = W~\text{almost everywhere}\}.
 \end{equation*}

 We equip $\widetilde{\Aut}(\Gamma)$ with the {\it weak topology} (e.g., \cite[weak topology]{hlm-erg}).
 \item Consider, next, the closed normal subgroup
 \begin{equation*}
 \widetilde{\Aut}(\Gamma)\trianglerighteq \widetilde{\Aut}_{\rm triv}(\Gamma)
 :=
\bigl\{g\in \widetilde{\Aut}(\Gamma)\mid W(g\, \cdot, \cdot) = W~\text{almost everywhere} \bigr\}.
 \end{equation*}
 Finally, set
 \begin{equation*}
 \Aut(\Gamma)
 :=
 \widetilde{\Aut}(\Gamma)/ \widetilde{\Aut}_{\rm triv}(\Gamma)
 \end{equation*}
 with the quotient topology.
 \end{itemize}
\end{Definition}

For pure graphons much of the quotienting (by trivial automorphisms, etc.) is not necessary.

\begin{Lemma}\label{le:purenotriv}
 If $\Gamma=(X,W,\mu)$ is a pure graphon, then $\widetilde{\Aut}_{\rm triv}(\Gamma)$ is trivial.
\end{Lemma}
\begin{proof}
 Let $\varphi\colon X\to X$ be a measure-preserving and a.e.\ $W$-preserving bijection with
 \begin{equation}\label{eq:wxwx}
 W(-,-) = W(\varphi -,-) \text{ almost everywhere}.
 \end{equation}
 We then have
 \begin{align*}
 \int_X r_W(x,\varphi x)\, \mathrm{d}\mu(x) &= \int_X\int_X |W(x,y)-W(\varphi x,y)|\, \mathrm{d}\mu(y) \mathrm{d}\mu(x)
 \quad\text{by Definition \ref{def:nbhd}}\\
 &= 0
 \quad\text{by \eqref{eq:wxwx}}.
 \end{align*}
 Since $r_W$ is assumed to be a metric and $\mu$ is fully supported, we have $x=\varphi x$ almost everywhere. In $\widetilde{\Aut}_{\rm triv}(\Gamma)$, we identify (by definition) automorphisms modulo measure zero, so we are done: $\varphi=\id$.
\end{proof}

Furthermore, the reason why $\Aut_{\rm LS}(\Gamma)$ functions so well is that it provides a section for the quotient by bijections that are trivial modulo measure zero.

\begin{Proposition}
 Let $\Gamma=(X,W,\mu)$ be a pure graphon and $\varphi\colon X\to X$ a measure-preserving bijection with
 \begin{equation*}
 W(\varphi -,\varphi - ) = W\ \text{almost everywhere}.
 \end{equation*}
 $\varphi$ is then congruent modulo zero to a unique element of $\Aut_{\rm LS}(\Gamma)$.
\end{Proposition}
\begin{proof}
 We are making two claims: existence and uniqueness. We will drop $\Gamma$ in the notation for automorphism groups: $\Aut_{\rm LS}=\Aut_{\rm LS}(\Gamma)$.

{\it Uniqueness.} The assertion to prove is that if two elements of $\Aut_{\rm LS}$ coincide modulo~0, then they coincide period. Or equivalently: if $\psi\in \Aut_{\rm LS}$ is trivial mod~0, then it is the identity.

 For $x\in X$, we have
 \begin{equation}\label{eq:autls}
 W(x,-) = W(\varphi x,\varphi-)\ \text{a.e.}
 \end{equation}
 because $\varphi\in \Aut_{\rm LS}$, and the assumption is that
 \begin{equation}\label{eq:phi1}
 \varphi = \id_X\ \text{a.e.}
 \end{equation}
 We then have
 \begin{align*}
 r_W(x,\varphi x) &= \int_X |W(x,y)-W(\varphi x,y)|\, \mathrm{d}\mu(y)
 \quad\text{by Definition~\ref{def:nbhd}}\\
 &= \int_X |W(x,y)-W(\varphi x,\varphi y)|\, \mathrm{d}\mu(y)
 \quad\text{by \eqref{eq:phi1}}\\
 &= \int_X |W(x,y)-W(x,y)|\, \mathrm{d}\mu(y)
 \quad\text{by \eqref{eq:autls}}\\
 &=0.
 \end{align*}
 Because $r_W$ is a metric, we must have $x=\varphi x$ and we are done: $x\in X$ was arbitrary.

{\it Existence.} The hypothesis implies that there is a full-measure (hence also dense) $Y\subseteq X$ such that $\varphi$ preserves $r_W$-distances between points in $Y$. But then, because $(X,r_W)$ is by assumption the completion of $(Y,r_W)$, the restriction $\varphi|_Y$ extends uniquely to an $r_W$ isometry $\varphi'$ on~$X$.

 We have now
 \begin{itemize}\itemsep=0pt
 \item modified $\varphi$ mod zero to $\varphi'$,
 \item so that the latter still preserves $W$ almost everywhere and is an $r_W$-isometry.
 \end{itemize}
 The conclusion that $\varphi'\in \Aut_{\rm LS}$ follows from \cite[Lemma~9]{ls} (or rather from its simpler analogue for $r_W$ rather than~$\overline{r}_W$).

This finishes the proof.
\end{proof}

In summary, we have the following.

\begin{Corollary}\label{cor:sameauto}
 For a pure graphon $\Gamma=(X,W,\mu)$, the automorphism groups $\Aut_{\rm LS}(\Gamma)$ and $\Aut(\Gamma)$ of Definitions {\rm \ref{def:origgrphn}} and {\rm \ref{def:autobis}} coincide. \qedhere
\end{Corollary}

As announced preceding Theorem \ref{th:prscr}, reasonable seemingly competing definitions of graphon automorphism groups thus yield the same objects.

\section{Rigidity}\label{se:rig}

\subsection{Graphon-rigid Lie groups}\label{subse:grrig-lie}

Recall \cite[Definition 30]{ls}'s notion of a \emph{Cayley graphon}: one whose underlying measure space is a compact second-countable (Hausdorff) group with its Haar probability measure, with
\begin{equation*}
 W(x,y)
 :=
 f\bigl(xy^{-1}\bigr)
 ,\qquad
 \bG\xrightarrow[\quad\text{measurable}\quad]{\quad f\quad}
 [0,1]
 ,\qquad
 f(x)=f\bigl(x^{-1}\bigr).
\end{equation*}
\cite[Theorem 37]{ls} notes that {\it weakly random} compact groups (i.e., those that have finitely many isomorphism classes of irreducible representations of dimension $\le d$ for every $d$) have the following ``permanence'' property.

\begin{Definition}\label{def:rig}
 A compact metrizable group $\bG$ is {\it graphon-rigid} if the limit of every convergent sequence of Cayley graphons on $\bG$ is again a Cayley graphon on $\bG$.
\end{Definition}

By contrast, \cite[Example 35]{ls} gives an example of Cayley graphons on the circle converging to a Cayley graphon on a 2-torus; essentially the same example is also discussed in \cite[text preceding Theorem 5]{sz-high}, and slightly adapted it reads as follows.

\begin{Example}\label{ex:t1t2}
 Let $\chi\colon \bS^1\to \bS^1$ be a generating character for the Pontryagin dual $\bZ\cong \widehat{\bS^1}$, and consider the Cayley graphons $\Gamma_n$ with
 \begin{itemize}\itemsep=0pt
 \item underlying Haar probability space $\bigl(\bS^1,\mu\bigr)$;
 \item and respective kernels
 \begin{equation*}
 W_n(x,y) = f_n\bigl(xy^{-1}\bigr),\qquad f_n:=\operatorname{Re}\left(\frac{2+\chi+\chi^n}4\right)
 \end{equation*}
 (`$\operatorname{Re}$' denoting the real part).
 \end{itemize}
 \cite[Example~35]{ls} observes that this sequence converges to the Cayley graphon on the torus
 \begin{equation*}
 \bT^2=\bS^1\times \bS^1
 \end{equation*}
 with kernel
 \begin{equation*}
 W(x,y)=f\bigl(xy^{-1}\bigr),\qquad f:=\operatorname{Re}\left(\frac{2+\chi_1+\chi_2}4\right),
 \end{equation*}
 where $\chi_i$, $i=1,2$, are the characters generating the Pontryagin duals of the two $\bS^1$ Cartesian factors. That agrees with \cite[discussion preceding Theorem 5]{sz-high}; we give one possible proof below, that will also help fit the example into a broader family.
\end{Example}

\begin{Remark}
 \cite[Example 35]{ls} actually uses the functions
 \begin{equation*}
 \operatorname{Im}\left(\frac{1+\chi+\chi^n}2\right)
 \quad\text{in place of}\quad
 \operatorname{Re}\left(\frac{2+\chi+\chi^n}4\right)
 \end{equation*}
 (i.e., {\it imaginary} rather than real parts), but this seems to not quite meet the requirements:
 \begin{itemize}\itemsep=0pt
 \item Per \cite[Definition~30]{ls}, one needs ``symmetric'' functions $f_n$ in the sense that $f_n\bigl(x^{-1}\bigr)=f_n(x)$. The real parts of the characters have this property, but their imaginary parts do not.
 \item And additionally, the kernels are meant to take values in $[0,1]$, whereas both the imaginary and the real part of $1+\chi+\chi^n$ will, in general, take some negative values.
 \end{itemize}
\end{Remark}

To unpack Example~\ref{ex:t1t2}, note first the following simple procedure for producing convergent graphon sequences; recall \cite[Section~8.1]{bog} that nets $(\mu_{\alpha})_{\alpha}$ of measures on topological spaces {\it converge weakly} to measures $\mu$ if
\begin{equation*}
 \int_X f\, \mathrm{d}\mu_{\alpha}\xrightarrow[\quad \alpha\quad]{} \int_X f\, \mathrm{d}\mu
\end{equation*}
for every bounded continuous function $f$ on $X$.

\begin{Lemma}\label{le:wconv}
 Let $X$ be a compact metrizable space and $W\colon X^2\to [0,1]$ a symmetric continuous function. If a sequence $\mu_n$ of probability measures converges weakly to $\mu$, then we correspondingly have graphon convergence
 \begin{equation*}
 \Gamma_n:=(X,W,\mu_n)\xrightarrow[\quad n\quad]{} (X,W,\mu)=:\Gamma.
 \end{equation*}
\end{Lemma}
\begin{proof}
 This is obvious from the characterization of convergence via morphism densities \cite[Section~2.1]{ls}:
 \begin{equation*}
 t(F,W_n)\xrightarrow[\quad n\quad]{} t(F,W)
 \end{equation*}
 for every finite simple graph $F=(V,E)$. Plainly, the hypothesis ensures that
 \begin{equation*}
 t(F,W_n):=\int_{X^V}\prod_{ij\in E}W(x_i,x_j)  \prod_{i\in V}\mathrm{d}\mu_n(x_i)
 \end{equation*}
 converges to the analogous quantity with $\mu$ in place of $\mu_n$.
\end{proof}

A slightly different take on Lemma~\ref{le:wconv} would be the following.

\begin{Lemma}\label{le:wconv-bis}
 Let
 \begin{equation*}
 \varphi_n\colon \ (X_n,\mu_n)\to (X,\mu)
 \end{equation*}
 be continuous measure-preserving functions between compact metrizable probability spaces and $(X,W,\mu)$ a graphon structure.

 Weak measure convergence
 \begin{equation*}
 \varphi_{n*}\mu_n\xrightarrow[\quad n \quad]{} \mu
 \end{equation*}
 then implies graphon convergence
 \begin{equation*}
 \bigl(X_n,W\circ\bigl(\varphi_n^{\times 2}\bigr),\mu_n\bigr)\xrightarrow[\quad n \quad]{} (X,W,\mu).
 \end{equation*}
\end{Lemma}
\begin{proof}
 Upon recalling \cite[Section~2.1]{ls} the weak isomorphism
 \begin{equation*}
 \bigl(X_n,W\circ\bigr(\varphi_n^{\times 2}\bigr),\mu_n\bigr) \cong (X,W,\varphi_{n*}\mu_n),
 \end{equation*}
 the conclusion follows from Lemma~\ref{le:wconv}.
\end{proof}

We can now place Example \ref{ex:t1t2} in broader context: it is an instance of Example \ref{ex:tt}, with $k=1$.

\begin{Example}\label{ex:tt}
 Let $\chi$ be a non-trivial character on $\bS^1$, as in Example \ref{ex:t1t2} (it need not be generating; any non-trivial character will do) and consider, for any
 $n_i\in \bZ_{>0}$, $1\le i\le k$,
 the function
 \begin{equation}\label{eq:kchars}
 f_{(n_i)}:=\operatorname{Re}\left(\frac 1{2(k+1)} \right) \bigl(k+1+\chi+\chi^{n_1}+\chi^{n_1 n_2}+\cdots+\chi^{n_1\cdots n_k}\bigr).
 \end{equation}
 Clearly, the $f_{(n_i)}$ are $[0,1]$-valued and symmetric, in the sense that they take equal values on $x\in \bS^1$ and $x^{-1}$. They thus induce kernels
$W_{(n_i)}(x,y):=f_{(n_i)}\bigl(xy^{-1}\bigr)$
 on~$\bS^1$. I claim that as the~$n_i$ all increase to infinity, the Haar-measure graphons
\smash{$
 \Gamma_{(n_i)}:=\bigl(\bS^1,W_{(n_i)},\mu\bigr)$}
 converge to the $(k+1)$-torus-based graphon $\Gamma=\bigl(\bT^{k+1},W,\mu_{k+1}\bigr)$, where $\mu_k$ is the Haar measure again and
$W(x,y) = f\bigl(xy^{-1}\bigr)$
 with $f$ being the symmetric function on $\bT^{k+1}$ constructed as in~\eqref{eq:kchars}
 \begin{equation*}
 f:=\operatorname{Re}\left(\frac 1{2(k+1)}\right)  (k+1+\chi_1+\cdots+\chi_{k+1} ),
 \end{equation*}
 the summands being $k+1$ non-trivial characters on the $k+1$ factors of the torus.

 To verify the claim, note first that it does not matter {\it which} \smash{$\chi_i\in \widehat{\bT^{k+1}}$} we select: all choices transform back to the one where all $\chi_i$ are generating via a component-wise measure preserving transformation of the torus, which will not alter the weak isomorphism class of the graphon. For that reason, we might as well take $\chi_i$ to be the $i^{\rm th}$-component copy of $\chi$.

 Next, observe that the functions $W_{(n_i)}$ are pullbacks of $W$ through the maps
 \begin{equation*}
 \bS^1\ni
 z
 \xmapsto[]{\quad \varphi_{(n_i)}\quad}
 \bigl(z,\, z^{n_1},\, z^{n_1n_2},\, \dots,\, z^{n_1\cdots n_k}\bigr)
 \in \bT^{k+1}.
 \end{equation*}
 It is not difficult to check now that the conditions of Lemma~\ref{le:wconv-bis} are met (except that we are working here with multi-sequences indexed by $(n_i)$ in place of a singly-indexed sequence): the images
\smash{$\varphi_{(n_i)}\bigl(\bS^1\bigr)\subset \bT^{k+1}$}
 converge to all of $\bT^{k+1}$ in the {\it Hausdorff metric} \cite[Definition~7.3.1]{bbi} on closed subsets of the latter, so any limit point of
$\varphi_{(n_i)*}\mu\in \mathrm{Prob}\bigl(\bT^{k+1}\bigr)$
 will be invariant under translations by the torus. The only such creature is the Haar measure~$\mu_{k+1}$, hence the conclusion (that the hypotheses of Lemma~\ref{le:wconv-bis} hold).
\end{Example}

\begin{Remark}\label{re:haustop}
 Defining the Hausdorff metric on the space $\mathcal{CL}(X)$ of closed subsets of a compact space~$X$ technically requires the ambient space to be metrized, but the metric will make no difference to the induced (compact \cite[Theorem~7.3.8]{bbi}) topology on $\mathcal{CL}(X)$.

 One way to see this is to observe that having fixed a metric $(X,d)$ inducing the underlying topology, convergence in $\mathcal{CL}(X)$ with respect to the induced Hausdorff metric can be defined in terms of just the {\it uniformity} induced by $d$ (\cite[Definition~7.1 and pp.~88--89]{jms-unif}). Because $X$ is compact Hausdorff that uniformity is in turn uniquely determined by the topology of $X$ \cite[Proposition~8.16]{jms-unif}.

 For an alternative take on the matter, note that the Hausdorff-metric topology on $\mathcal{CL}(X)$ (associated to any metric topologizing $X$) can also be recovered intrinsically, with no reference to a metric, as the topology introduced by Fell on \cite[p.~472]{fell-top} (see also, e.g., \cite[Chapter~II, Section~2]{am} for a recollection): a subbasis of open sets consists of the sets
 \begin{equation*}
 \{Y\subseteq X\text{ closed}\mid Y\cap K=\varnothing,\, Y\cap U\ne\varnothing,\, \forall U\in \cF\}.
 \end{equation*}
 as $K\subseteq X$ ranges over compact subsets and $\cF$ over finite collections of open subsets of $X$.

 The definition makes sense for arbitrary topological spaces $X$, and makes $\mathcal{CL}(X)$ compact Hausdorff whenever~$X$ is {\it locally} compact Hausdorff \cite[Theorem~1]{fell-top}. This is not the case, in general, for Hausdorff-metric topologies, but the various topologies {\it do} agree when $X$ is compact (rather than just locally so); this being a simple enough exercise, we omit the proof.

 In any event, when $X$ is furthermore a compact {\it group}, the collection $\mathcal{CLG}(X)$ of closed sub{\it groups} is closed and hence also compact \cite[p.~474, (IV)]{fell-top}. In fact, \cite[Theorem~1]{fell-top} makes it reasonable to work with $\mathcal{CL}(X)$ (or $\mathcal{CLG}$, when the ambient space is a group) even when $X$ is compact but {\it not} metrizable.
\end{Remark}

We need a few auxiliary observations for later.

\begin{Lemma}\label{le:biinv}
 Let $\bG$ be a compact second-countable group. There is a continuous symmetric function $f\colon\bG\to [0,1]$ such that
$\bG\ni x\mapsto f(x\, \cdot)$
 is injective.
\end{Lemma}
\begin{proof}
 The topology of $\bG$ is induced by a metric $d$, bi-invariant in the sense that~\cite[Corollary~A4.19]{hm4}
 \begin{equation*}
 d(gx,gy) = d(x,y) = d(xg,yg),\qquad \forall g,x,y\in \bG.
 \end{equation*}
To conclude, set
$f(x):=d(x,1)$.
 For $x\ne y\in \bG$, the translates $f(x\ \cdot)$ and $f(y\ \cdot)$ take the value~0 at~$x^{-1}$ and~$y^{-1}$, respectively, and nowhere else.
\end{proof}

And a consequence thereof is the following.

\begin{Lemma}\label{le:nondeg}
 Every second-countable compact group has a Cayley-graphon structure meeting the requirements of Proposition~{\rm \ref{pr:cpct}}.
\end{Lemma}
\begin{proof}
 Indeed, construct one using a function $f$ as in Lemma~\ref{le:biinv}.
\end{proof}

It will be useful to distill a general principle operative in Example~\ref{ex:tt}.

\begin{Lemma}\label{le:hausconv}
 Let $\Gamma:=(\bG,W,\mu_{\bG})$ be a Cayley graphon on a compact metrizable group and
$ \varphi_n\colon  \bG_n\to \bG$
 morphisms of compact metrizable groups.
 If $\varphi_n(\bG_n)$ converge to $\bG$ in the Hausdorff distance induced by any metric topologizing $\bG$, then we have convergence
 \begin{equation*}
 \Gamma_n:=\bigl(\bG_n,W\circ \varphi_n^{\times 2},\mu_{\bG_n}\bigr)\xrightarrow[\quad n\quad]{} \Gamma
 \end{equation*}
 of Cayley graphons.
\end{Lemma}
\begin{proof}
 As in Example \ref{ex:tt}, any limit-point of the sequence of pushforward probabilities $\varphi_{*n}\mu_{\bG_n}$ will be translation-invariant under the Hausdorff limit
$\bG = \lim_n \varphi_n (\bG_n)$,
 so must be the Haar measure $\mu_{\bG}$. We can then apply Lemma~\ref{le:wconv-bis} to conclude.
\end{proof}

Incidentally, we can supplement Example \ref{ex:tt} with non-Lie analogues.

\begin{Example}\label{ex:zp}
 Let $p$ be a prime number, and denote by
 \begin{equation}\label{eq:invlim}
 \bZ_p := \varprojlim_n \bZ/p^n
 \end{equation}
 the profinite group of {\it $p$-adic integers} \cite[Example 1.28]{hm4}.

 Equip the circle $\bS^1$ with a Cayley-graphon structure $\bigl(\bS^1,W,\mu\bigr)$ as in Lemma~\ref{le:nondeg}, and take for the maps $\varphi_n$ of Lemma~\ref{le:hausconv} the morphisms
\smash{$\bZ_p\twoheadrightarrow \bZ/p^n\hookrightarrow \bS^1$}
 obtained by first surjecting along the canonical map provided by~\eqref{eq:invlim} and then embedding the finite cyclic group as the \smash{$(p^n)^{\rm th}$} roots of unity.

 The hypothesis of Lemma~\ref{le:hausconv} clearly holds, and hence so does its conclusion. We thus have $\bZ_p$-Cayley graphons converging to a circle-based Cayley graphon. The latter is not $\bZ_p$-Cayley though: the $p$-adics would have to act transitively (and continuously, Lemma~\ref{le:actcont}) on the manifold $\bS^1$, which would entail that the action factor through a~Lie \cite[Proposition~1.6.5]{tao-hilb} and hence finite quotient. Naturally, finite groups do not act transitively on the circle.
\end{Example}

To make sense of Theorem \ref{th:lierig}, recall \cite[Definition~9.5]{hm4} that a compact connected group~$\bG$ is {\it semisimple} if it equals its {\it commutator subgroup}~$\bG'$; this is the subgroup generated by commutators
$[x,y]:=xyx^{-1}y^{-1}$, $x,y\in \bG$.
According to \cite[Theorem 9.2]{hm4}, it is also just the {\it set} of such commutators, and is always connected and closed in $\bG$.

We also remind the reader that {\it simple} compact groups \cite[Definition~9.88]{hm4} are those having no proper, non-trivial normal subgroups (or equivalently \cite[Theorem~9.90]{hm4}, no such {\it closed} subgroups).

For compact {\it Lie} groups, it turns out that \cite[Theorem~37]{ls} in fact provides a characterization of graphon-rigidity. The following result collects that along with a number of other such characterizations.

\begin{Theorem}\label{th:lierig}
 Let $\bG$ be a compact Lie group and $\bG_0$ its connected identity component. The following conditions are equivalent:
 \begin{enumerate}[$(a)$]\itemsep=0pt
 \item\label{item:5} $\bG$ is graphon-rigid.
 \item\label{item:6} $\bG_0$ is graphon-rigid.
 \item\label{item:7} $\bG_0$ is semisimple.
 \item\label{item:8} $\bG_0$ has finite center.
 \item\label{item:9} $\bG_0$ does not admit a surjection $\bG\to \bS^1$ onto the circle group.
 \item\label{item:10} Modulo a finite central subgroup, $\bG_0$ is a product of finitely many compact connected simple Lie groups.
 \item\label{item:13} More precisely, there are compact connected simple Lie groups $\bS_i$, $1\le i\le n$, finite central subgroups $\bN_i\le \bS_i$, and surjections
 \begin{equation}\label{eq:snd}
 \prod_{i=1}^n \bS_i
 \xrightarrow[\quad]{}
 \bG_0
 \xrightarrow[\quad]{}
 \prod_{i=1}^n \bS_i/\bN_i.
 \end{equation}
 \item\label{item:14} $\bG$ is weakly random.
 \item\label{item:15} $\bG_0$ is weakly random.
 \end{enumerate}
\end{Theorem}

We first dispose of the relationship between conditions \ref{item:5} and \ref{item:6}, in slightly more general form than needed. Wherever compact groups act transitively on graphons $(X,W,\mu)$, we will assume the latter pure; it follows from \cite[Corollary 16]{ls} (and its proof) that they are also compact, with topology induced by either of the metrics $r_W$ or $\overline{r}_W$.

\begin{Proposition}\label{pr:finix}
 Let $\bG$ be a compact metrizable group whose connected component $\bG_0\le \bG$ has finite index. If $\bG$ is graphon-rigid, so is $\bG_0$.
\end{Proposition}
\begin{proof}
 Let
 \begin{equation*}
 \Gamma_n \xrightarrow[\quad n\quad ]{} \Gamma,\quad \Gamma_n\text{ $\bG_0$-Cayley}
 \end{equation*}
 be a convergent sequence of graphons.

 As all $\Gamma_n:=(X_n,W_n,\mu_n)$ come equipped with transitive (measure-preserving) $\bG_0$-actions, we can form the {\it induced} probability $\bG$-spaces in the sense of \cite[Appendix G]{kech-erg}:
 \begin{itemize}\itemsep=0pt
 \item Set
 \begin{equation*}
 \overline{X_n}:=\mathrm{Ind}_{\bG_0}^{\bG}X_n := X_n\times \bG/\bG_0
 \end{equation*}
 as a space.
 \item The probability measure is
 \begin{equation*}
 \overline{\mu_n}:=\mu_n\times(\text{normalized counting measure on }\bG/\bG_0).
 \end{equation*}
 \item The kernel $\overline{W_n}$ is a copy of the original $W_n$ on each individual copy of $X_n$. For points
 \begin{equation*}
 (x,s)\in X_n\times\{s\},\qquad (x',s')\in X_n\times\{s'\},\qquad s\ne s'\in \bG/\bG_0,
 \end{equation*}
 we can set, for instance,
 \begin{equation*}
 \overline{W_n}((x,s),(x',s'))=0;
 \end{equation*}
 this will do (making the resulting induced graphon pure) provided $X_n$ was not a singleton with $W\equiv 0$. We assume $X_n$ and $X$ are at any rate not singletons, since otherwise there would be (almost) nothing to prove.
 \item The (also transitive) $\bG$-action is as explained in \cite[Appendix G]{kech-erg}: having chosen a finite set $S\subset \bG$ of representatives for the cosets $\bG/\bG_0$, the action of $g\in \bG$ is
 \begin{equation*}
 g(x,s) = (hx,s')\qquad \text{for the unique }gs=s'h,\ h\in \bG_0,\ s,s'\in S.
 \end{equation*}
 \end{itemize}
 The analogous construction applies to $\Gamma:=(X,W,\mu)$ to produce
 \begin{equation*}
 \overline{\Gamma}:=\mathrm{Ind}_{\bG_0}^{\bG}\Gamma = \bigl(\overline{X},\overline{W},\overline{\mu}\bigr),
 \end{equation*}
 and it is easy to see that this induction procedure is continuous:
 \begin{equation*}
 \overline{\Gamma_n} \xrightarrow[\quad n\quad ]{} \overline{\Gamma}.
 \end{equation*}
 The hypothesis says that the limit $\overline{\Gamma}$ is $\bG$-Cayley. The orbits of the {\it connected} compact group~$\bG_0$ are each contained within single ``slices''
 \begin{equation}\label{eq:xs}
 X\cong X\times \{s\},\qquad s\in S\cong \bG/\bG_0,
 \end{equation}
 so each orbit has measure $\le \frac 1{[\bG:\bG_0]}$. But the orbits also cover $X$ (by $\bG$-transitivity), and there must be at least $[\bG:\bG_0]$ of them. It follows that each $\bG_0$-orbit constitutes a single leaf~\eqref{eq:xs}; as those are isomorphic to the original $\Gamma$, we are done.
\end{proof}

\begin{proof}[Proof of Theorem \ref{th:lierig}]
 We do this in stages.

Part 1: \ref{item:7} through \ref{item:13} are mutually equivalent. This is standard compact-group structure theory: according to \cite[Theorem~9.24]{hm4}, say, an arbitrary compact connected $\bG_0$ can be sandwiched as in~\eqref{eq:snd}, with finite central kernels, except that perhaps there are additional torus factors on the two sides (apart from the semisimple factors). Those torus factors are absent precisely when $\bG_0$ is semisimple \cite[Theorem~9.19 and Corollary~9.20]{hm4}.

Part 2: \ref{item:14} $\Leftrightarrow$ \ref{item:15}. The identity component $\bG_0\le \bG$ has finite index because $\bG$ is compact and Lie. The usual ``Mackey machine'' applies \cite[Theorems~4.64 and~4.65]{kt} to conclude that
 \begin{itemize}\itemsep=0pt
 \item every irreducible $\bG$-representation is a sum of at most $[\bG:\bG_0]$ conjugates under $\bG$ of an irreducible $\bG_0$-representation;
 \item and every irreducible $\bG_0$-representation arises in this fashion, as a summand of an irreducible $\bG$-representation.
 \end{itemize}
 This gives a relation between irreducible $\bG$- and $\bG_0$-representations: two such are related precisely when one contains the other. Related representations have dimensions differing by a factor of at most $[\bG:\bG_0]$, and the relation is at most $[\bG:\bG_0]$-to-one.

 It follows, then, that a cap of $N\in \bZ_{>0}$ on the number of irreducible $\bG$-representations of dimension $\le d$ will give a cap of $[\bG:\bG_0]N$ on the number of irreducible $\bG_0$-representations of dimension $\le \frac d{[\bG:\bG_0]}$ and vice versa.

Part 3: \ref{item:5} $\Rightarrow$ \ref{item:6}. This follows from Proposition~\ref{pr:finix}, since for compact {\it Lie} groups the identity component is also open and hence has finite index.

 This, so far, breaks up the various conditions into three aggregates:
 \begin{enumerate}[(1)]\itemsep=0pt
 \item\label{item:16} \ref{item:5} $\Rightarrow$ \ref{item:6},
 \item\label{item:17} \ref{item:7} up to \ref{item:13}, all mutually equivalent,
 \item\label{item:18} \ref{item:14} and \ref{item:15}, all mutually equivalent.
 \end{enumerate}
 To establish {\it inter}-group equivalence, note first that~\ref{item:18} implies~\ref{item:16} by \cite[Theorem~37]{ls}, so only two other implications are needed.

Part 4: \ref{item:17} $\Rightarrow$ \ref{item:18}. This says that semisimple compact connected Lie groups are weakly random; the claim follows easily, say, from familiar Lie-group/algebra representation theory via, say, the {\it Weyl dimension formula} for representations \cite[Theorem~IX.6.1]{sim-rep}.

Part 5: \ref{item:16} $\Rightarrow$ \ref{item:17}. Or the contrapositive: if $\bG_0$ {\it does} admit a surjection onto $\bS^1$, then it cannot be graphon-rigid.

 By \cite[Theorem~9.24]{hm4}, $\bG_0$ surjects with finite central kernel onto a Lie group of the form
 \begin{equation}\label{eq:tss}
 \bT^{\ell}\times\prod_{i=1}^n \bS_i,\qquad \ell\ge 1,\quad \bS_i\text{ simple}.
 \end{equation}
 That finite kernel will make no material difference to the argument, so we will simplify matters and assume~$\bG_0$ is~\eqref{eq:tss} to begin with.

 Now consider a somewhat larger group
$
 \bH:=\bT^{\ell+1}\times\prod_{i=1}^n \bS_i$
 (note the higher torus dimension), and equip it with a Cayley graphon structure $\Gamma:=(\bH,W,\mu)$ as in Lemma~\ref{le:nondeg}, so that the conditions of Proposition~\ref{pr:cpct} hold.

 Next, consider maps
$\varphi_n\colon \bG_0\to \bH$
 that
 \begin{itemize}\itemsep=0pt
 \item operate as the identity on the $\bS_i$ and all but one of the $\ell$ circle factors in $\bT^{\ell}$;
 \item and wind one of the circle factors ever more tightly around $\bT^2$ via, say,
 \begin{equation*}
 \bS^1\ni z\xmapsto[]{\quad \varphi_n\quad} (z,z^n)\in \bT^2
 \end{equation*}
 \end{itemize}
 (as in Example~\ref{ex:tt}). It follows from Lemma~\ref{le:hausconv} that the Cayley graphons
$
 \Gamma_n:=\bigl(\bG_0,\allowbreak {W\circ \varphi_n^{\times 2}},\mu_{\bG_0}\bigr)$, $n\in \bZ_{>0}$,
 converge to $\Gamma$. I claim, however, that $\bG_0$ cannot act transitively on their limit $\Gamma$.

 Indeed, such an action would be continuous on the compact metric space $(\bH,r_W)$ (Lem\-ma~\ref{le:actcont}), so would identify $\bH$ homeomorphically with the homogeneous space
$\bG_0/\bG_{0,p}$ (for some $p\in \bH$)
 by the (closed, hence Lie \cite[Theorem 20.12]{lee}) stabilizer group $\bG_{0,p}\le \bG_0$ of a point $p\in \bH$. But that homogeneous space is a topological manifold of dimension \cite[Theorem~21.17]{lee}
$
 \dim \bG_0-\dim \bG_{0,p}\le \dim \bG_0$,
while $\bH$ is a manifold of strictly larger dimension $\dim \bG_0+1$. That the two cannot be homeomorphic then follows from dimension invariance \cite[Theorem~17.26]{lee}.

 This concludes the proof.
\end{proof}

It is perhaps worth noting at this stage that for general (non-Lie) compact metrizable groups graphon-rigidity does {\it not} imply weak randomness: the latter is thus the strictly stronger condition.\looseness=-1

\begin{Example}\label{ex:z2s}
 Let $\bG = (\bZ/2)^{\aleph_0}$, the product of countably infinitely many copies of $\bZ/2$. Naturally, $\bG$ is not weakly random: being abelian, {\it all} of its (infinitely many) irreducible representations are 1-dimensional. I claim that it is nevertheless graphon-rigid.

 To see this, note that by the argument employed in the proof of \cite[Theorem 39]{ls}, any limit graphon of a sequence of $\bG$-Cayley graphons will be acted upon transitively by an abelian group~$\bA$ of exponent~$\le 2$ (i.e., such that $x^2=1$ for all elements~$x$).

 Any compact group of exponent $\le 2$ is a product of copies of $\bZ/2$, and since all groups in sight are metrizable, at most countably many copies at that. $\bG$ will then surject onto~$\bA$, hence the conclusion.
\end{Example}

Recall the following useful device from \cite[Section~4.2]{ls}.

\begin{Definition}\label{def:trunc}
 Let $\Gamma:=(X,W,\mu)$ be a graphon and $r > 0$. The {\it truncation} $[\Gamma]_{r}$ is the graphon obtained from $\Gamma$ as follows:
 \begin{itemize}\itemsep=0pt
 \item Regard $W$ as a {\it Hilbert--Schmidt operator} \cite[Section~4]{halm-sund} on $L^2:=L^2(X,\mu)$ by
 \begin{equation*}
 (Wf)(x):=\int_X W(x,y)f(y)\, \mathrm{d}\mu(y),\qquad f\in L^2.
 \end{equation*}
 \item Let $\lambda$ be an eigenvalue of $W$. Denoting by $P_{\lambda}$ the orthogonal projection onto the (finite-dimensional) $\lambda$-eigenspace of $W$, the {\it truncation} $[W]_r$ is
$
 \sum_{|\lambda|\ge r}\lambda P_{\lambda}$.
 It is again a \emph{kernel} in the sense of \cite[Section~7.1]{lvs-bk}, i.e., a bounded symmetric measurable function (the boundedness being noted in \cite[Section~4.1]{ls}); it is not, however, necessarily non-negative (see Remark \ref{res:post.trunc}\,\ref{item:res:post.trunc:not.grphn}).
 \item Then set
 $[\Gamma]_r:=(X,[W]_r,\mu)$.
 \end{itemize}
\end{Definition}

\begin{remarks}\label{res:post.trunc}\quad
 \begin{enumerate}[(1)]\itemsep=0pt
 \item\label{item:res:post.trunc:not.grphn} Note that the $[W]_r$ just introduced need not be \emph{graphon} structures (though they are referred to as such in \cite[Section~4.2]{ls}), for there is no reason, in principle, why the truncated kernels should be non-negative; I am grateful to one of the anonymous referees for this caution.

 Numerical exploration confirms this: one can produce random non-negative-entry symmetric $3\times 3$ matrices and spectrally truncate away the smallest eigenvalue so as to produce some negative entries in the result. Given the possible negativity of $[W]_r$, these are what one might call \emph{generalized} graphons: based on arbitrary kernels rather than $[0,1]$-valued such.\looseness=-1

 The same notion of convergence does make sense, as do automorphism groups. The latter, moreover, easily transport between the plain-graphon and generalized-graphon settings: a~bounded symmetric $W$ can always be shifted into positivity and scaled into $[0,1]$, retaining the automorphism group.

 \item\label{item:res:post.trunc:fin.approx} The truncation procedure of Definition \ref{def:trunc} gives a kind of canonical ``finite-dimensional approximation'' of the graphon $\Gamma$. For $r>0$, the purification of the truncated $[\Gamma]_r$ can be realized as a subspace of
$
 V_r:=\bigoplus_{|\lambda|\ge r}(\text{$\lambda$-eigenspace of $W$})$,
 equipped with the kernel
$
 V_r^2\ni (x,y)\mapsto \sum_i \lambda_i x_iy_i$
 for an enumeration $(\lambda_i)$ of the relevant eigenvalues $|\lambda|\ge r$ and choice of coordinates
$
 V_r\ni x = (x_i)_i\in \bR^{\dim V_r}$.

 Furthermore, it follows from \cite[Lemma 17]{ls} that if $r>0$ is ``generic enough'' (specifically, not an eigenvalue of $\Gamma$), then for any convergent sequence $\Gamma_n\to \Gamma$ the truncations $[\Gamma_n]_r$ also converge to $[\Gamma]_r$, the measures associated with a subsequence thereof converge weakly to that attached to $[\Gamma]_r$, etc.
 \end{enumerate}
\end{remarks}

Example \ref{ex:z2s} replicates for {\it connected} compact groups as well.

\begin{Example}\label{ex:qhats}
 Consider the Pontryagin dual $\widehat{\bQ}$ of the additive group of rationals. Since $\bQ$ is a~filtered colimit (in the category of discrete abelian groups) of copies of $\bZ$, we have a description%
 \begin{equation}\label{eq:qhatts}
 \widehat{\bQ}\cong \varprojlim \bS^1
 \end{equation}
 as a limit (in the category of compact abelian groups), where the connecting maps are $z\mapsto z^n$ for various $n$; the target group will be an infinite product $\bG=\widehat{\bQ}^{\aleph_0}$. Once more, being abelian, it has only 1-dimensional irreducible representations. We argue that it too is graphon-rigid however, much along the same lines as in Example \ref{ex:z2s}, adopting \cite[Theorem 39]{ls}'s proof strategy.

 Let $\Gamma:=(X,W,\mu)$ be the limit of a sequence of $\bG$-Cayley graphons $\Gamma_n$. For fixed $\alpha>0$ assume, as in the proof of \cite[Theorem 39]{ls}, that the truncations $[\Gamma_n]_{\alpha}$ as well as $[\Gamma]_\alpha$ are all supported in the same Euclidean space $\bR^d$. The $[\Gamma_n]_{\alpha}$ are acted upon transitively by $\bG$, hence by images of $\bG$ via morphisms $\bG\to \bO(d)$; such images must be
 \begin{itemize}\itemsep=0pt
 \item Lie because they are compact subgroups of $\bO(d)$;
 \item and abelian connected, being quotients of $\bG$ (which has both properties).
 \end{itemize}
 In short, $[\Gamma_n]_{\alpha}$ admit transitive actions by tori $\bT_n$. Passing to a subsequence if necessary, we can assume convergence
$ \lim_n \bT_n =: \bT_{\alpha}=\bT\le \text{orthogonal group }\bO(d)$
 in the Hausdorff topology of Remark~\ref{re:haustop} on the compact space $\mathcal{CLG}(\bO(d))$ of closed subgroups of $\bO(d)$. The transitivity of the $\bT_n$ on $[\Gamma_n]_{\alpha}$ ensures that of $\bT$ on $[\Gamma]_{\alpha}$, and I next claim that $\bT$ too must be a torus (equivalently, given $\bT$'s being a Lie group, connected and abelian). To see this, note that
 \begin{itemize}\itemsep=0pt
 \item any two arbitrary elements of $\bT$ are arbitrarily approximable by pairs of elements in the $\bT_n$s, so must commute;
 \item and, were $\bT$ a \emph{dis}connected closed subgroup of $\bO(d)$, a sufficiently small neighborhood of any non-identity connected component will avoid all $\bT_n$ with sufficiently large $n$, all $\bT_n$ being connected.
 \end{itemize}
 Now let $\alpha\to 0$, so that the inverse limit $\bH$ of the $\bT_{\alpha}$ acts transitively on the original graphon $\Gamma$ (cf.\ \cite[proof of Theorem~39, last paragraph]{ls}). The desired conclusion, that a limit of $\bG$-Cayley graphons admits a transitive $\bG$-action, will have been achieved once we observe that $\bH$ admits surjections $\bG\xrightarrowdbl{}\bH$. Indeed: being a countable colimit of free, finitely generated abelian groups, the Pontryagin dual $\widehat{\bH}$ will be a countable torsion-free discrete abelian group. All such embed into
$\bQ^{\oplus \aleph_0}\cong \widehat{\bG}$
 (as an application of \cite[Lemma 57.5 and Theorem 57.13]{cr}, say), hence the claimed existence of surjective morphisms $\bG\xrightarrowdbl{}\bH$.
\end{Example}

A familiar phenomenon, when studying symmetries of objects parametrized by a space, is that of the {\it upper semicontinuity} of their groups of automorphisms. A case in point is that of Riemannian structures: consider a compact smooth manifold~$M$ and the space $\cat{Riem}(M)$ of Riemannian structures, appropriately topologized \cite[Section~1]{ebin}. The diffeomorphism group~$\bD$ of $M$ acts on $\cat{Riem}(M)$, and the isotropy group
\begin{equation*}
 \bD_{(M,g)} := \{\gamma\in \bD\mid \gamma\text{ fixes }(M,g)\in \cat{Riem}(M)\}
\end{equation*}
of a point (i.e., Riemannian structure) is nothing but the isometry group $\mathrm{Aut}(M,g)$: always compact and Lie. It turns out that for points $(M,g')$ close to $(M,g)$ the automorphism groups $\mathrm{Aut}(M,g')$ are ``almost contained'' in $\mathrm{Aut}(M,g)$ \cite[Theorem 8.1]{ebin}: they are conjugate in $\bD$ to subgroups of $\mathrm{Aut}(M,g)$ by elements close to the identity. This upper semicontinuity phenomenon does {\it not} obtain for graphons (at least not in anything like the same fashion).

\begin{Example}\label{ex:smallsph}
 Consider {\it any} compact probability metric space $(X,d,\mu)$ as a graphon, and scale the distance towards 0. This produces a net converging to a singleton equipped with the zero kernel; all members $(X,cd,\mu)$ of the net (for scalars $c>0$) have the same automorphism group, which may well be non-trivial, while naturally, the limit has trivial automorphism group.
\end{Example}

Or more interestingly (albeit branching out into generalized graphons) is the following.

\begin{Example}\label{ex:qhat}
 According to Theorem \ref{th:prscr}, there is a graphon $\Gamma:=(X,W,\mu)$ having $\widehat{\bQ}$ as an automorphism group
$\Aut(\Gamma)\cong \widehat{\bQ}$.
 The eigenvalue-truncations of~$\Gamma$ (as in \cite[Sections~4.2 and~4.4]{ls} and Example \ref{ex:qhats}) all have compact Lie automorphism groups containing quotients of~$\widehat{\bQ}$ acting non-trivially. The only compact Lie quotients of $\widehat{\bQ}$ are also connected and abelian, so~they must be tori \cite[Chapter~II, Exercise~C.2]{helg}. But they cannot be {\it higher-dimensional} tori~$\bT^d$, $d\ge 2$, because
\smash{$\widehat{\bT^d}\cong \bZ^d$}, $d\ge 2$,
 does not embed into \smash{$\bQ\cong \widehat{\widehat{\bQ}}$} \cite[Theorem~7.63]{hm4}.

 In short: $\Gamma$, with automorphism group $\widehat{\bQ}$, is a limit of generalized graphons with automorphism groups containing copies of $\bS^1$. It is enough to note now that the latter cannot embed into \smash{$\widehat{\bQ}$}, because its dual $\bZ$ is not a quotient of $\bQ$.
\end{Example}

\subsection{Stable compact-group images}\label{subse:stabim}

By way of motivation for the subsequent material, we have the following.

\begin{Remark}\label{re:tbad2ways}
 Starting with a $\widehat{\bQ}$-Cayley graphon, its truncations are $\bS^1$-Cayley (arguing as in Example~\ref{ex:qhat}: they admit transitive actions by Lie quotients of $\widehat{\bQ}$, etc.). It follows, then, that the circle $\bS^1$ rigid-deficient through two qualitatively distinct mechanisms:
 \begin{enumerate}[(a)]\itemsep=0pt
 \item\label{item:wrap} it can ``wrap around'' higher-dimensional tori increasingly densely, as in Example~\ref{ex:tt} (and \cite[Example~35]{ls});
 \item\label{item:surj} and it self-surjects
$ \bS^1\ni z\mapsto z^n\in \bS^1$
 with increasingly large kernel, giving rise to the inverse limit \eqref{eq:qhatts}.
 \end{enumerate}
\end{Remark}

The following notion is intended to isolate and address the type of behavior noted in item~\ref{item:wrap} of Remark~\ref{re:tbad2ways}.

\begin{Definition}\label{def:downrig}
 A compact group $\bG$ is {\it image-rigid} if for every compact {\it Lie} group $\bH$ the images of morphisms $\bG\to \bH$ constitute a closed set of groups in the space $\mathcal{CLG}(\bH)$ of Remark \ref{re:haustop}.
\end{Definition}

An equivalent formulation for image-rigidity, occasionally useful, is the following.

\begin{Lemma}\label{le:denseall}
 A compact group $\bG$ is image-rigid in the sense of Definition~{\rm \ref{def:downrig}} if and only if it satisfies the following property:
 $\bG$ surjects onto every compact Lie group $\bH$ belonging to the closure in $\mathcal{CLG}(\bH)$ of images of morphisms $\bG\to \bH$.
\end{Lemma}
\begin{proof}
 That Definition \ref{def:downrig} is formally stronger than the present requirement is clear, so only the other implication is interesting.

 Assume then that $\bG$ has the property in the statement, and consider morphisms
 \begin{equation*}
 \varphi_n\colon \ \bG\to \bK,\qquad \bK\text{ compact Lie}
 \end{equation*}
 with
 \begin{equation*}
 \varphi_n(\bG)\xrightarrow[\quad n\quad]{} \bH\le \bK\qquad\text{in}\quad \mathcal{CLG}(\bK).
 \end{equation*}
 By \cite[Theorem~1 and its Corollary]{mz-nearby}, for sufficiently large $n$ we can find
 \begin{equation*}
 \bK\ni x_n\xrightarrow[\quad n\quad]{}1,\qquad x_n \varphi_n(\bG) x_n^{-1}\le \bH.
 \end{equation*}
 Replacing $\varphi_n$ with their respective conjugates $x_n \cdot \varphi(-)\cdot x_n^{-1}$, we can place ourselves entirely within $\bH$; the property in the statement then applies to deliver the conclusion (that $\bG$ surjects onto $\bH$).
\end{proof}

Examples will be familiar by now.

\begin{Example}
 As noted repeatedly, a circle $\bS^1$ is not image-rigid: the limit in $\mathcal{CLG}(\bT^2)$ of the images of
$\bS^1\ni z\mapsto (z,z^n)\in \bT^2$ as $n\to \infty$
 is the full 2-dimensional torus, but of course {\it it} is not an image of the circle.
 By the same token, none of the tori $\bT^n$ are image-rigid: they all wrap around larger-dimensional tori in the same fashion.
\end{Example}

In fact, for Lie groups we have just recovered another facet of Theorem \ref{th:lierig}.

\begin{Proposition}\label{pr:lieredux}
 For a compact Lie group $\bG$ the following conditions are equivalent $($and hence also equivalent to those of Theorem~{\rm \ref{th:lierig})}:
 \begin{enumerate}[$(a)$]\itemsep=0pt
 \item\label{item:4} $\bG$ has semisimple connected component $\bG_0$.
 \item\label{item:11} For every compact Lie group $\bH$ the space of morphisms $\bG\to \bH$ is compact in the compact-open topology.
 \item\label{item:12} $\bG$ is image-rigid.
 \end{enumerate}
\end{Proposition}

Before moving on to the proof, we pause to note that things are again somewhat more complicated in the non-Lie case.

\begin{Example}\label{ex:inf.tor}
 An {\it infinite}-dimensional torus $\bT^{\aleph_0}$ is image-rigid: its images $\bT^{\aleph_0}\to \bH$ in a~compact Lie group $\bH$ are all compact and connected and hence (finite-dimensional) tori \cite[Theo\-rem~4.2.4]{de}, as are limits thereof in $\mathcal{CLG}(\bH)$; $\bT^{\aleph_0}$ surjects onto every finite-dimensional torus, hence the claim.

 On the other hand, $\bT^{\aleph_0}$ is not graphon-rigid. We noted in Example \ref{ex:qhat} and Remark \ref{re:tbad2ways} that $\widehat{\bQ}$-Cayley graphons will be limits of $\bS^1$-Cayley ones, which in turn are also $\bT^{\aleph_0}$-Cayley. On the other hand, the infinite torus cannot act continuously and transitively on $\widehat{\bQ}$.

 One way to verify this last claim is to observe that $\bT^{\aleph_0}$ is path-connected, while $\widehat{\bQ}$ is not. These claims, in turn, follow from \cite[Theorem~8.62]{hm4}: according to that result the quotient of a compact abelian group $\bG$ by the (possibly non-closed) path component of the identity is the {\it ext group} \cite[Definition A1.51]{hm4} $\mathrm{Ext}\bigl(\widehat{\bG},\bZ\bigr)$. Now,
 \begin{itemize}\itemsep=0pt
 \item For $\bG=\bT^{\aleph_0}$, we have
$\mathrm{Ext}\bigl(\widehat{\bG},\bZ\bigr)\cong \mathrm{Ext}\bigl(\bZ^{\oplus \aleph_0},\bZ\bigr) = \{0\}$
 (e.g., by \cite[Proposition~A1.14]{hm4}), hence path connectedness.
 \item On the other hand, for $\bG=\widehat{\bQ}$ we have
$\mathrm{Ext}\bigl(\widehat{\bG},\bZ\bigr)\cong \mathrm{Ext}(\bQ,\bZ) \cong \bR\ne \{0\}$
 by \cite{wie-qz}.
 \end{itemize}
\end{Example}

Some preparation will help simplify the argument.

\begin{Lemma}\label{le:conncomp}
 For any compact group $\bG$ the set of closed subgroups with abelian connected component is closed in $\mathcal{CLG}(\bG)$.
\end{Lemma}
\begin{proof}
 Expressing an arbitrary compact group as a limit of compact {\it Lie} groups \cite[Corollary~2.43]{hm4}, it is enough to focus on these: throughout the proof, then, $\bG$ will be Lie.

 Every closed subgroup $\bH\le \bG$ with abelian identity component is a product
$\bH = \bH_0 \cdot \bF$, $\bF\le \bH$ is finite,
 by \cite[Theorem~6.10]{hm4}, so a limit of a net $\bH_n$ can be recovered (perhaps after passing to a subnet) as the product between the (connected, abelian) limit of~$\bH_{n,0}$ and that of finite subgroups $\bF_n\le \bH_n$. It thus suffices to prove the main claim for {\it finite} groups.

 That, in turn, follows from the fact that $\bG$ being Lie and hence a matrix group \cite[Corollary~9.58]{hm4}, there is a positive integer $N$ such that every finite subgroup thereof has an abelian subgroup of index $\le N$ \cite[Theorem~36.13]{cr}. The limit of a convergent sequence of finite groups~$\bF_n$ will again have an abelian subgroup of index $\le N$, and we are done.
\end{proof}

The lemma in turn helps with the following result, intended to reduce (part of) the problem to groups with abelian identity component. As in Theorem \ref{th:lierig} (and the discussion preceding it), denote by $\bG_0$ the identity component of a compact group $\bG$ and by $\bG'\le \bG$ the (algebraic) commutator subgroup; it will often be closed automatically, e.g., when $\bG$ is connected \cite[Theorem~9.2]{hm4} or Lie \cite[Theorem~6.11]{hm4}.

\begin{Proposition}\label{pr:abcom}
 If a compact group $\bG$ is image-rigid in the sense of Definition~{\rm \ref{def:downrig}}, so are $\bG/\bG_0'$ and $\bG_{ab}:=\bG/\bG'$.
\end{Proposition}
\begin{proof}
 Write $\overline{\bG}$ for either $\bG/\bG_0'$ or $\bG_{ab}$, and consider morphisms $\varphi_n\colon \overline{\bG}\to \bH$ into a compact Lie group with respective images
$
 \bH_n:=\varphi_n\bigl(\overline{\bG}\bigr)$.
 A~limit $\bH_{\infty}$ of $\bH_n$ in $\mathcal{CLG}(\bH)$ is an image of $\bG$ by assumption, and we need to argue that it is in fact an image of~$\overline{\bG}$.

 When $\overline{\bG}=\bG_{ab}$ this is obvious: every $\bH_n$ is then abelian, so $\bH_{\infty}$ will again be so; a morphism from $\bG$ onto it will then factor through the abelianization $\bG_{ab}$. We thus focus on the case $\overline{\bG}=\bG/\bG_0'$.

 Since $\varphi_n$ maps $\bG$ onto $\bH_n$, it also maps~$\bG_0$ {\it onto} $\bH_{n,0}$ \cite[Lemma~9.18]{hm4}. Because it annihilates~$\bG_0'$ by assumption, $\bH_n$ have the property that their connected components are abelian. This must hold of $\bH_{\infty}$ too by Lemma~\ref{le:conncomp}, and we conclude as before: a surjective morphism $\bG\to \bH_{\infty}$ restricts on $\bG_0$ to a morphism through the latter's abelianization
$\bG_{0,ab}:=\bG_0/\bG_0'$,
 so must annihilate $\bG_0'$.
\end{proof}

\begin{proof}[Proof of Proposition~\ref{pr:lieredux}]
 That \ref{item:11} is formally stronger than \ref{item:12} is obvious, so we focus on the other two implications.

\ref{item:12} $\Rightarrow$ \ref{item:4}. Suppose $\bG_0$ is not semisimple; we will argue that the quotient $\bG/\bG_0'$ by the derived subgroup of its identity component is not image-rigid, so neither is $\bG$ by Proposition~\ref{pr:abcom}.

 The abelianization $\bG_0/\bG_0'$ is a non-trivial ($d$-dimensional, say) torus $\bT^d$, $d\ge 1$; to simplify the notation, we assume $\bG_0'$ vanishes to begin with and $\bG_0=\bT^d$. The finite quotient $\bF:=\bG/\bT^d$ acts on $\bT^d$ by conjugation, and the extension
$\{1\}\to \bT^d\to \bG\to \bF\to \{1\}$
 can be regarded as an element $\alpha$ of the cohomology group $H^2\bigl(\bF,\bT^d\bigr)$ \cite[Section~I.1]{moore-ext}. Now make~$\bF$ act on
$ \bT^{2d}\cong \bT^d\times \bT^d$
 diagonally, by doubling up the original action, and consider the $\bF$-equivariant morphism
 \begin{equation*}
 \bT^d\ni z\xmapsto[]{\quad\varphi_n\quad} (z,z^n)\in \bT^{2d}.
 \end{equation*}
 It pushes $\alpha$ forward to an element of $H^2\bigl(\bF,\bT^{2d}\bigr)$, which in turn corresponds to an extension $\bG_1$ of $\bF$ by $\bT^{2d}$ instead:
 \begin{equation*}
 \begin{tikzpicture}[auto,baseline=(current bounding box.center)]
 \path[anchor=base]
 (0,0) node (ll) {$\{1\}$}
 +(2,.6) node (lu) {$\bT^d$}
 +(2,-.6) node (ld) {$\bT^{2d}$}
 +(4,.6) node (mu) {$\bG$}
 +(4,-.6) node (md) {$\bG_1$}
 +(6,0) node (r) {$\bF$}
 +(7,0) node (rr) {$\{1\}$.}
 ;
 \draw[->] (ll) to[bend left=6] node[pos=.5,auto] {$\scriptstyle $} (lu);
 \draw[->] (ll) to[bend right=6] node[pos=.5,auto] {$\scriptstyle $} (ld);
 \draw[->] (lu) to[bend left=0] node[pos=.5,auto] {$\scriptstyle $} (mu);
 \draw[->] (ld) to[bend left=0] node[pos=.5,auto] {$\scriptstyle $} (md);
 \draw[->] (mu) to[bend left=6] node[pos=.5,auto] {$\scriptstyle $} (r);
 \draw[->] (md) to[bend right=6] node[pos=.5,auto] {$\scriptstyle $} (r);
 \draw[->] (r) to[bend right=0] node[pos=.5,auto] {$\scriptstyle $} (rr);
 \draw[->] (lu) to[bend left=0] node[pos=.5,auto] {$\scriptstyle \varphi_n$} (ld);
 \draw[->] (mu) to[bend left=0] node[pos=.5,auto] {$\scriptstyle $} (md);
 \end{tikzpicture}
 \end{equation*}
 Clearly, as $n\to\infty$ the images of $\bG$ will converge to the higher-dimensional $\bG_1$, which cannot itself be an image of $\bG$. This proves the claim that $\bG$ is not image-rigid.

\ref{item:4} $\Rightarrow$ \ref{item:11}. Embedding $\bH$ into a unitary group $\bU(n)$ \cite[Corollary 9.58]{hm4}, it will be enough to consider morphisms $\bG\to \bU(n)$ or, equivalently, $n$-dimensional unitary $\bG$-representations.

 Since we saw in Theorem \ref{th:lierig} that the semisimplicity of $\bG$ entails weak randomness, there are finitely many irreducible $\bG$-representations in each dimension $\le n$. It follows that the space of morphisms $\bG\to \bU(n)$ decomposes as a disjoint union of finitely many orbits under $\bU(n)$-conjugation, finishing the proof.
\end{proof}

Compact Lie groups are very close to connected: their connected components are open and hence cofinite. At the other extreme stand {\it profinite} groups: inverse limits of finite groups or, equivalently, totally disconnected compact groups \cite[Theorem 1.34]{hm4}. For abelian ones, at least, the two notions of rigidity converge (together with a complete characterization of those that enjoy that property).

\begin{Theorem}\label{th:profab}
 For a profinite abelian group $\bG$, the following conditions are equivalent:
 \begin{enumerate}[$(a)$]\itemsep=0pt
 \item\label{item:21} $\bG$ is image-rigid.
 \item\label{item:22} We have
 \begin{equation}\label{eq:prodzns}
 \bG\cong \prod_{i\in I}\bZ/n_i,\qquad (n_i)_i\text{ bounded}.
 \end{equation}
 \item\label{item:23} $\bG$ is graphon-rigid.
 \end{enumerate}
\end{Theorem}
\begin{proof}
 We complete an implication cycle.

\ref{item:21} $\Rightarrow$ \ref{item:22}. Having an isomorphism \eqref{eq:prodzns} is equivalent to the discrete Pontryagin dual $\widehat{\bG}$ being of the form
 \begin{equation}\label{eq:ghatbdd}
 \widehat{\bG}\cong \bigoplus_i \bZ/n_i,\qquad (n_i)_i\text{ bounded}.
 \end{equation}
 The dual $\widehat{\bG}$ is in any event torsion (being dual to a profinite group \cite[Corollary 8.5]{hm4}), and \eqref{eq:ghatbdd} is in turn equivalent (by \cite[Theorem 6]{kap}, say) to $\widehat{\bG}$ being of {\it bounded order}: there is an upper bound on the order of an element of $\widehat{\bG}$.

 Suppose, now, that $\widehat{\bG}$ is {\it not} of bounded order. This means that we can find characters
$\chi\colon \bG\to \bS^1$
 of arbitrarily large order, and hence with arbitrarily large image. But then images of such characters will converge in $\mathcal{CLG}(\bS^1)$ to the whole circle, which of course is not an image of the profinite group~$\bG$ (images of compact profinite groups are profinite \cite[Exercise~E1.13]{hm4}). Image rigidity is thus violated.

\ref{item:22} $\Rightarrow$ \ref{item:23}. This is a slight elaboration on (and adaptation of) Example~\ref{ex:z2s}: a limit of Cayley graphons over \eqref{eq:prodzns} will be acted upon transitively by a countable inverse limit of quotients of that group; I claim that the product in question surjects onto any such metrizable inverse (this suffices, graphon automorphism groups being metrizable). The rest of the present argument is devoted to confirming that claim.

 First, it will simplify matters (notationally if nothing else) to specialize the $n_i$ somewhat: quotients and limits alike decompose as products indexed by (finitely many, by boundedness) prime numbers $p$ of \emph{$p$-primary components} (i.e., \cite[p.~5]{kap} quotients of $\bG$ whose elements all have $p$-power orders). We will thus assume without loss of generality that $\bG\cong \prod\bZ/p^{k_j}$ for a~bounded exponent set $\{k_j\}$.

For an arbitrary $p$-primary product $\bP$ of cyclic $p$-groups, write $\sharp_m \bP$ for the cardinality of the set of cyclic factors surjecting onto $\bZ/p^m$. The inverse limit $\bH\cong \varprojlim_s \bH_s$ of metric quotients $\bG\xrightarrowdbl{} \bH_s$ we want to argue $\bG$ surjects onto is precisely such a product (as is $\bG$), and the surjection will obtain provided $\sharp_m \bG\ge\sharp_m \bH$ for all $m$. That this is so follows from the countability of the limit and the fact that $\sharp_m \bG\ge \sharp_m \bH_s$ for all~$s$.

\ref{item:23} $\Rightarrow$ \ref{item:21}. The claim goes through without the abelianness assumption; we relegate it to Proposition~\ref{pr:prof2rig}.
\end{proof}

\begin{Proposition}\label{pr:prof2rig}
 A compact graphon-rigid group $\bG$ is image-rigid if it is either
 \begin{enumerate}[$(a)$]\itemsep=0pt
 \item\label{item:24} profinite;
 \item\label{item:25} or connected.
 \end{enumerate}
\end{Proposition}
\begin{proof}
 In evaluating image rigidity, Lemma~\ref{le:denseall} provides a compact Lie group $\bH$ and morphisms
$
 \varphi_n\colon   \bG\to \bH$,
 whose images converge to all of $\bH$ in $\mathcal{CLG}(\bH)$. The argument now bifurcates.

\ref{item:24}: profinite groups. Failing image rigidity, Lemma~\ref{le:conncomp} allows us to assume that $\bH$ has abelian non-trivial identity component $\bH_0$. But then the argument of Example~\ref{ex:zp} applies to produce $\bG$-Cayley graphons converging to a $\bH$-Cayley graphon, while at the same time showing that the latter does not admit a transitive $\bG$-action.

 \ref{item:25}: connected groups. This time $\bH$ must be connected. We consider two cases.

(I) $\bH$ is non-abelian. It is then the isometry group of a left-invariant metric $(\bH,d)$ by \cite[Theorem~1.3]{niem} (which requires non-abelianness as one of its stipulations). The Haar measure~$\mu$ will then complete a graphon structure
$\Gamma:=(\bH,d,\mu)$
 (perhaps after rescaling $d$ so that it takes values in the unit interval), and the morphisms $\varphi_n$ give us a sequence of $\bG$-Cayley graphons converging to $\Gamma$ as in Lemma~\ref{le:hausconv}. Graphon rigidity means that we have a morphism
$\bG\to \Aut(\Gamma)\cong \bH$.
 Its image must act transitively on $\bH$ by translation, so the morphism in question is onto; hence the conclusion.

(II) $\bH$ is abelian. So that it must be a torus
$\bT^{\ell}\cong \bR^{\ell}/\bZ^{\ell}$, $\ell\in \bZ_{>0}$
  (being compact, connected, abelian and Lie \cite[Theorem~4.2.4]{de}). Now equip $\bT^{\ell}$ with its metric $d$ induced by the standard Euclidean norm on $\bR^{\ell}$ and its Haar measure~$\mu$. It is easy to see that the $\bT^{\ell}$-translations constitute precisely the identity component $\Aut(\Gamma)_0$ of
$\Aut(\Gamma)$, $\Gamma:=\bigl(\bT^{\ell},d,\mu\bigr)$.
 Now conclude as previously: the morphisms $\varphi_n$ generate $\bG$-Cayley graphon structures converging to $\Gamma$, so that graphon rigidity will give a morphism
$ \varphi\colon \bG\to \Aut(\Gamma)$
 whose image
 \begin{itemize}\itemsep=0pt
 \item on the one hand is contained in $\Aut(\Gamma)_0\cong \bT^{\ell}$ because $\bG$ is connected;
 \item and on the other hand, operates transitively on the torus by translation.
 \end{itemize}
 In short, $\varphi$ is onto.
\end{proof}

\section{Graphing automorphisms}\label{se:graphing}

There is some indication\footnote{Lov{\'a}sz L.,
 Personal communication, Email to the author, 21~May 2022.}
 that a treatment of \emph{graphing} \cite[Section~2]{lvs-grphng} automorphism groups along parallel lines might not be without interest. We record a few observations on that theme in the present section for future reference, starting with some recollections.

\begin{Definition}\label{def:phing}
 A {\it graphing} $\Gamma=(V,\cB,E,\mu)$ consists of a standard probability space $(V,\cB,\mu)$ and a graph $(V,E)$ with no loops such that
 \begin{itemize}\itemsep=0pt
 \item $E\subset V\times V$ is a Borel set;
 \item the degree of $(E,V)$ is uniformly bounded by some positive integer $D$;
 \item and the degree-symmetry condition
 \begin{equation*}
 \int_A\deg_B(x)\, \mathrm{d}\mu(x) = \int_B\deg_A(x)\, \mathrm{d}\mu(x)
 \end{equation*}
 holds for all Borel sets $A,B\in \cA$, where $\deg_B(x)$ denotes the number of elements in $B$ connected to $x$.
 \end{itemize}
\end{Definition}

Following \cite[Section~2]{lvs-grphng}, the {\it edge measure} of a graphing $\Gamma=(V,\cB,E,\mu)$ is defined by
\begin{equation*}
 \eta(A\times B)=\eta_{\Gamma}(A\times B):=\int_A\deg_B(x)\, \mathrm{d}\mu(x)
\end{equation*}
for $A,B\in \cB$.

In light of Remark \ref{re:mod0} (specifically, part \ref{item:3}), the most straightforward route to defining the automorphism group of a graphing is via the Borel measure algebra.

\begin{Definition}\label{def:phing-auto}
 Let $\Gamma=(V,\cB,E,\mu)$ be a graphing. Its {\it automorphism group} $\Aut(\Gamma)$ consists of those automorphisms of the measure algebra $(\cB,\mu)$ that also preserve the edge measure $\eta_{\Gamma}$ on~${\cB\times \cB}$.
\end{Definition}

By contrast to the automorphism group of a graphon, which is compact by \cite[Theorem~10]{ls}, that of a graphing need not be.

\begin{Example}\label{ex:ncpct}
 If the edge-set $E$ is empty, the automorphism group $\Aut(V,\cB,E,\mu)$ is nothing but the full automorphism group of a the standard probability space $(V,\cB,\mu)$. If that space is non-atomic (and hence can be identified with the unit interval with its Lebesgue measure), or even as soon as its non-atomic component is non-empty, that group is certainly non-compact.
\end{Example}

Recall the following \cite[Section~2]{lvs-grphng}.

\begin{Definition}
 A {\it full subgraphing} of a graphing $(V,\cB,E,\mu)$ is a full subgraph supported by a~full-measure Borel subset of~$V$.
\end{Definition}

It is immediate that full subgraphings retain the automorphism groups of their larger ambient graphings.

\begin{Proposition}
 For a full subgraphing $\Gamma'$ of a graphing $\Gamma=(V,\cB,E,\mu)$, we have an isomorphism $\Aut(\Gamma')\cong \Aut(\Gamma)$. \qedhere
\end{Proposition}

As explained in \cite[Section~2]{lvs-grphng}, full subgraphings are both {\it locally-globally equivalent} and plain {\it locally equivalent} to the original graphings in the sense of \cite[Definition~3.3]{hls}. Local equivalence alone does not suffice to preserve the automorphism group.

\begin{Example}\label{ex:loc.eq.alone}
 The graph(ing) consisting of a single vertex and no edges is locally equivalent to the two-vertex discrete graph assigning probability $\frac 12$ to each of the vertices (by \cite[Example~3.5]{hls}, for instance). Clearly though, the first graph has trivial automorphism group whereas in the second we can interchange the vertices (measure-preservingly).
\end{Example}

On the other hand, recall \cite[Proposition 7.11]{hls} that graphings are locally equivalent precisely when they are \emph{bi-locally isomorphic}: there is a graphing admitting local isomorphisms into both (local isomorphisms being possibly non-invertible measure-preserving maps between vertex sets, restricting to isomorphisms between connected components of the respective graphings).

\subsection*{Acknowledgments}

This work is partially supported by NSF grant DMS-2001128.
I am grateful for helpful comments and pointers to the literature from L.~Lov\'asz, as well as the anonymous referees' careful comments and suggestions.

\pdfbookmark[1]{References}{ref}
\LastPageEnding

\end{document}